\documentclass[1p]{elsarticle}
\usepackage{amssymb}
\usepackage{amsfonts}

\newtheorem{theorem}{Theorem}
\newtheorem{conjecture}[theorem]{Conjecture}

\newtheorem{lemma}[theorem]{Lemma}

\newtheorem{remark}{Remark}

\newproof{pf}{Proof}

\begin{document}
\title{Distant sum distinguishing index of graphs with bounded minimum degree}

\author{Jakub Przyby{\l}o\fnref{fn1,fn2}}
\ead{jakubprz@agh.edu.pl, phone: 048-12-617-46-38,  fax: 048-12-617-31-65}

\fntext[fn1]{Financed within the program of the Polish Minister of Science and Higher Education
named ``Iuventus Plus'' in years 2015-2017, project no. IP2014 038873.}
\fntext[fn2]{Partly supported by the Polish Ministry of Science and Higher Education.}

\address{AGH University of Science and Technology, al. A. Mickiewicza 30, 30-059 Krakow, Poland}

\begin{abstract}
For any graph $G=(V,E)$ with maximum degree $\Delta$ and without isolated edges, and a positive integer $r$,
by $\chi'_{\Sigma,r}(G)$ we denote the $r$-distant sum distinguishing index of $G$.
This is the least integer $k$ for which a proper edge colouring $c:E\to\{1,2,\ldots,k\}$ exists
such that $\sum_{e\ni u}c(e)\neq \sum_{e\ni v}c(e)$ for every pair of distinct vertices $u,v$ at distance at most $r$ in $G$.
It was conjectured that $\chi'_{\Sigma,r}(G)\leq (1+o(1))\Delta^{r-1}$ for every $r\geq 3$.
Thus far it has been in particular proved that $\chi'_{\Sigma,r}(G)\leq 6\Delta^{r-1}$ if $r\geq 4$.
Combining probabilistic and constructive approach, we show that this can be improved to $\chi'_{\Sigma,r}(G)\leq (4+o(1))\Delta^{r-1}$
if the minimum degree of $G$ equals at least $\ln^8\Delta$.
\end{abstract}


\maketitle

\section{Introduction}

Integer edge colourings were initiated in the paper of Chartrend et al.~\cite{Chartrand},
where the graph invariant \emph{irregularity strength}, $s(G)$, was introduced as a possible measure of the `level of irregularity' of a graph $G$.
This referred to the well known phenomenon in graph theory that there are no \emph{irregular graphs}, understood as graphs whose all vertices have pairwise distinct degrees (see also~\cite{ChartrandErdosOellermann} for possible alternative definitions of irregularity in graphs), except the trivial 1-vertex case. For a given graph $G=(V,E)$, $s(G)$ is defined as the least $k$ for which
one is able to construct an irregular multigraph (defined analogously as in the case of graphs above) of $G$ by multiplying some of its edges -- each at most $k$ times.
In terms of integer colourings, the same value is equivalently defined as
the least $k$ so that an edge colouring $c:E\to\{1,2,\ldots,k\}$ exists attributing every vertex $v\in V$ a distinct \emph{weighted degree} defined as:
$$d_c(v):=\sum_{e\ni v}c(e).$$
This we shall also call the \emph{sum at} $v$, see e.g.
\cite{Aigner,Bohman_Kravitz,Lazebnik,Dinitz,Faudree2,Faudree,Frieze,KalKarPf,Lehel,MajerskiPrzybylo2,Nierhoff,Przybylo,irreg_str2}
for a few out of a vastness of results concerning $s(G)$, which also gave rise to a whole discipline devoted to investigating this and other related problems.
One of the most intriguing direct descendants of the irregularity strength is its local correspondent, 
where we necessarily require an inequality $d_c(u)\neq d_c(v)$ to hold only for adjacent vertices $u,v$ in $G$.
The least $k$ admitting a colouring $c:E\to\{1,2,\ldots,k\}$ with such a feature we shall denote by $s_1(G)$.
In the first paper~\cite{123KLT} concerning this the authors conjectured that $k=3$ suffices for every connected graph of order at least $3$.
This presumption is commonly referred to as the \emph{1--2--3 Conjecture} nowadays.
This was investigated e.g. in~\cite{Louigi30,Louigi,123with13}.
The best thus far general result is however the upper bound
$s_1(G)\leq 5$ from~\cite{KalKarPf_123}.
A generalization of this concept, forming a link between $s_1(G)$ and $s(G)$, was introduced in~\cite{Przybylo_distant}.
Let $d(u,v)$ denote the distance of vertices $u,v$ in $G$.
We shall call $u$ and $v$, $r$-neighbours if $1\leq d(u,v)\leq r$ in $G$, where $r$ is a positive integer.
For every vertex $v$ in $G$, the set of its $r$-neighbours shall be denoted by $N^r(v)$, and we set $d^r(v)=|N^r(v)|$.
The least $k$
so that an edge colouring $c:E\to\{1,2,\ldots,k\}$ exists with $d_c(u)\neq d_c(v)$ for every $r$-neighbours $u,v\in V$ in $G$
is denoted by $s_r(G)$
(note it would be justified to set $s_{\infty}(G)=s(G)$ in the same spirit),
see e.g.~\cite{Przybylo_distant} and~\cite{Przybylo_distant_edge_probabil} for a few results concerning this concept,
which refers to the known distant chromatic numbers (see~\cite{DistChrSurvey} for a survey of this topic in turn).

In this paper we shall investigate a related problem referring to distant chromatic numbers.
Given a positive integer $r$ and a graph $G=(V,E)$ without isolated edges,
the \emph{$r$-distant sum distinguishing index} of $G$,
denoted by $\chi'_{\Sigma,r}(G)$, is the least integer $k$ such that there exists a \emph{proper} edge colouring $c:E\to\{1,2,\ldots,k\}$
which \emph{sum-distinguishes} $r$-neighbours in $G$, i.e. such that
$d_c(u)\neq d_c(v)$ for every $u,v\in V$ with $1\leq d(u,v)\leq r$.
In~\cite{Przybylo_distant_edge_proper} the following conjecture,
approximating the investigated lower bounds discussed e.g. in~\cite{Przybylo_distant,Przybylo_distant_edge_proper},
was posed.
\begin{conjecture}[\cite{Przybylo_distant_edge_proper}]\label{Conjecture1+o(1)}
For every integer $r\geq 3$ and each graph $G$ without isolated edges of maximum degree $\Delta$,
$\chi'_{\Sigma,r}(G)\leq (1+o(1))\Delta^{r-1}$. 
\end{conjecture}
It was also conjectured under the same conditions, that $\chi'_{\Sigma,2}(G)\leq (2+o(1))\Delta$~\cite{Przybylo_distant_edge_proper},
and that $\chi'_{\Sigma}(G)=\chi'_{\Sigma,1}(G)\leq \Delta+2$ for every connected graph $G$ of order at least $3$ non-isomorphic to $C_5$~\cite{FlandrinMPSW}.
Thus far for $r\geq 4$, the following is known.
\begin{theorem}[\cite{Przybylo_distant_edge_proper}]\label{mainTheoremDistPropSum}
Let $G$ be a graph without isolated edges and with maximum degree
$\Delta\geq 2$, and let $r\geq 4$. Then
$\chi'_{\Sigma,r}(G)\leq 6\Delta^{r-1}$.
\end{theorem}
Upper bounds of orders conjectured above are also known for $r=2,3$, but with slightly worse multiplicative constants than in Theorem~\ref{mainTheoremDistPropSum} above, see~\cite{Przybylo_distant_edge_proper},
while the upper bound of the form $\chi'_\Sigma(G)\leq (1+o(1))\Delta(G)$ was proved in~\cite{Przybylo_asym_optim} and~\cite{Przybylo_asymptotic_note},
see also~\cite{BonamyPrzybylo,DongWang_mad,FlandrinMPSW,Przybylo_CN_1,Przybylo_CN_2,WangChenWang_planar} for other results concerning the case $r=1$.
In this paper we combine probabilistic approach with a special constructive algorithm in order to provide the following improvements of the best known upper bounds for all $r\geq 4$ from Theorem~\ref{mainTheoremDistPropSum}, under assumption that the minimum degree of a graph is larger than some poly-logarithmic function of the maximum degree. (The value of this function, which seems unavoidable within our approach, could still be optimized -- we did not try to do this for the sake of clarity of the presentation).

\begin{theorem}\label{przybylo_main_th-sum-proper-probabil}
For every integer $r\geq 4$ there exists a constant $\Delta_0$ such that for each
graph $G$ with maximum degree $\Delta\geq \Delta_0$ and minimum degree $\delta\geq \ln^8\Delta$,
$$\chi'_{\Sigma,r}(G) < 4\Delta^{r-1}\left(1+ \frac{3}{2\ln\Delta}\right) 
+384,$$
hence $\chi'_{\Sigma,r}(G)\leq (4+o(1))\Delta^{r-1}$ for all graphs with $\delta\geq \ln^8\Delta$ and without isolated edges.
\end{theorem}

\section{Probabilistic Tools and Preliminary Lemmas}
The following standard tools of the probabilistic method shall be applied: 
the Lov\'asz Local Lemma, see e.g.~\cite{AlonSpencer},
the Chernoff Bound, see e.g.~\cite{JansonLuczakRucinski}
(Th. 2.1, page 26)
and Talagrand's Inequality, see e.g.~\cite{MolloyReed_GoodTalagrand}.
Details follow.
\begin{theorem}[\textbf{The Local Lemma}]
\label{LLL-symmetric}
Let $A_1,A_2,\ldots,A_n$ be events in an arbitrary pro\-ba\-bi\-li\-ty space.
Suppose that each event $A_i$ is mutually independent of a set of all the other
events $A_j$ but at most $D$, and that ${\rm \emph{\textbf{Pr}}}(A_i)\leq p$ for all $1\leq i \leq n$. If
$$ ep(D+1) \leq 1,$$
then $ {\rm \emph{\textbf{Pr}}}\left(\bigcap_{i=1}^n\overline{A_i}\right)>0$.
\end{theorem}
\begin{theorem}[\textbf{Chernoff Bound}]\label{ChernofBoundTh}
For any $0\leq t\leq np$,
$${\rm\emph{\textbf{Pr}}}({\rm BIN}(n,p)>np+t)<e^{-\frac{t^2}{3np}}~~{and}~~{\rm\emph{\textbf{Pr}}}({\rm BIN}(n,p)<np-t)<e^{-\frac{t^2}{2np}}\leq e^{-\frac{t^2}{3np}}$$
where ${\rm BIN}(n,p)$ is the sum of $n$ independent Bernoulli variables, each equal to $1$ with probability $p$ and $0$ otherwise.
\end{theorem}

\begin{theorem}[\textbf{Talagrand's Inequality}]\label{TalagrandsInequalityTotal}
Let $X$ be a non-negative random variable 
determined by
$l$ independent trials $T_1,\ldots,T_l$. 
Suppose there exist constants $c,k>0$ such that for every set of possible outcomes of the trials, we
have:
\begin{itemize}
\item[1.] changing the outcome of any one trial can affect $X$ by at most $c$, and
\item[2.] for 
each $s>0$, if $X\geq s$ then there is a set of at most $ks$ trials whose outcomes certify that $X\geq s$. 
\end{itemize}
Then for any 
$t\geq 0$, we have
\begin{equation}\label{TalagrandsInequality}
{\rm\emph{\textbf{Pr}}}(|X-{\rm\emph{\textbf{E}}}(X)|>t+20c\sqrt{k{\rm\emph{\textbf{E}}}(X)}+64c^2k)\leq 4e^{-\frac{t^2}{8c^2k(\mathbf{E}(X)+t)}}.
\end{equation}
\end{theorem}

We note that
knowing that 
$\mathbf{E}(X)\leq h$ 
we may also 
apply Talagrand's Inequality e.g. to
the variable $Y=X+h-\mathbf{E}(X)$, with $\mathbf{E}(Y)=h$ to obtain the following 
counterpart of~(\ref{TalagrandsInequality})
provided that the assumptions of Theorem~\ref{TalagrandsInequalityTotal}
hold for $X$: 
$$\mathbf{Pr}(X>h+t+20c\sqrt{kh}+64c^2k) \leq \mathbf{Pr}(Y>h+t+20c\sqrt{kh}+64c^2k) \leq 4e^{-\frac{t^2}{8c^2k(h+t)}}.$$
Similarly, the Chernoff Bound can be applied e.g. when we know that
$X$ is a sum of $n\leq k$ 
random independent Bernoulli variables, each equal to $1$ with probability at most $q$, 
to prove that $\mathbf{Pr}(X>kq+t)< e^{-\frac{t^2}{3kq}}$ (if $t\leq\lfloor k\rfloor q$).

In order to prove our main result we shall need the following observation.

\begin{lemma}\label{SparseSpanningSubgraphInC}
If $\Delta$ is large enough, then
for every graph $G'$ of maximum degree $\Delta'\leq \Delta$ and with minimum degree $\delta'\geq \frac{1}{2}\ln^5\Delta$,
there exists a spanning subgraph $F'$ of $G'$ with $d_{F'}(v)\leq \frac{d_{G'}(v)}{\ln^3\Delta}$ for each $v\in V(G')$.
\end{lemma}

\begin{pf}
We assume that $\Delta$ is large enough so that all inequalities within the proof below hold.
Independently for every vertex $v\in V(G')$ choose one of its incident edges, each with equal probability,
and denote the subgraph induced in $G'$ by the set of all the chosen edges by $F'$.
We shall show that with positive probability such $F'$ complies with our requirements.
For every $v\in V(G')$ denote by $X_v$ the random variable representing the number
of all edges incident with $v$ and chosen to $E(F')$ by any of the neighbours of $v$ in $G'$,
and note that $d_{F'}(v)\leq X_v+1$ (as at most one more edge incident with $v$ in $G'$ might be chosen to $E(F')$ by $v$ itself).
Note that for any given vertex $v\in V(G')$ and its neighbour $u\in N_{G'}(v)$,
the probability that $uv$ was chosen by $u$ equals 
$\frac{1}{d_{G'}(u)}\leq \frac{2}{\ln^5\Delta}$,
hence 
$${\mathbf E}(X_v)\leq \frac{2d_{G'}(v)}{\ln^5\Delta}\leq \frac{d_{G'}(v)}{2\ln^3\Delta}-\frac{1}{2}.$$
By the Chernoff Bound (with $t=\frac{d_{G'}(v)}{2\ln^3\Delta}-\frac{1}{2} \geq \frac{\ln^2\Delta}{5}$) we thus obtain that
\begin{equation}\label{Chernoff_ineq_subgraph}
{\mathbf Pr} (X_v > \frac{d_{G'}(v)}{\ln^3\Delta} - 1) < e^{-\frac{\ln^2\Delta}{15}} < \frac{1}{\Delta^3}.
\end{equation}
As any event $X_v > \frac{d_{G'}(v)}{\ln^3\Delta} - 1$ is mutually independent of all other events $X_{v'} > \frac{d_{G'}(v')}{\ln^3\Delta} - 1$
with $d(v,v')>2$, 
i.e. all except at most $(\Delta')^2\leq \Delta^2$,
by~(\ref{Chernoff_ineq_subgraph}) and the Lov\'asz Local Lemma we may conclude that with positive probability
for every $v\in V(G')$, $X_v\leq \frac{d_{G'}(v)}{\ln^3\Delta} - 1$, hence $d_{F'}(v) \leq \frac{d_{G'}(v)}{\ln^3\Delta}$.
A desired $F'$ must thus exist. \qed
\end{pf}

We shall also need to guarantee a special ordering of the vertices of a graph $G=(V,E)$.
For any linear ordering of $V$ and a vertex $v\in V$, a neighbour or $r$-neighbour 
of $v$ which precedes it in the
ordering shall be called a \emph{backward neighbour} or \emph{$r$-neighbour}, resp., of $v$.
The remaining ones in turn shall be referred to as \emph{forward neighbours} or \emph{$r$-neighbours}, resp., of $v$,
while the edges joining $v$ with its forward or backward neighbours shall be called \emph{forward} or \emph{backward}, resp., as well.
For any subset $S\subset V$, let also $N_-(v)$, $N_-^r(v)$, $N_S(v)$, $N_S^r(v)$ denote the sets of all backward neighbours, backward $r$-neighbours, neighbours in $S$ and $r$-neighbours in $S$ of $v$, respectively.
Set finally
$d_-(v)=|N_-(v)|$, $d_-^r(v)=|N_-^r(v)|$, $d_S(v)=|N_S(v)|$, $d_S^r(v)=|N_S^r(v)|$, and for any subset of edges $E_0\subseteq E$, $d_{E_0}(v)=|\{u\in N(v):uv\in E_0\}|$.\\

The following lemma was proved in~\cite{Przybylo_distant_edge_probabil}.
Here we only outline the main ideas behind its proof -- 
the remaining part of the argument
can however be reconstructed by an interested reader,
as in general it is based on a similar combination of the Chernoff Bound and Local Lemma as the (less complex) proof of Lemma~\ref{SparseSpanningSubgraphInC} above.

\begin{lemma}[\cite{Przybylo_distant_edge_probabil}]\label{MainSequencingLemma}
There exists a constant $\Delta'_0$ such that for every graph $G=(V,E)$ with maximum degree $\Delta\geq \Delta'_0$ and minimum degree $\delta\geq \ln^8\Delta$,
there is an assignment attributing every vertex $v\in V$ a distinct real number in $X_v\in[0,1]$ such that if we denote:
\begin{eqnarray}
A&=&\left\{v:X_v<\frac{1}{\ln^2\Delta}\right\},\nonumber\\
B&=&\left\{v:\frac{1}{\ln^2\Delta}\leq X_v\leq 1-\frac{1}{\ln^{3}\Delta}\right\},\nonumber\\
C&=&\left\{v:X_v> 1-\frac{1}{\ln^{3}\Delta}\right\}\nonumber
\end{eqnarray}
and order the vertices in $V$ into the sequence $v_1,v_2,\ldots,v_n$ consistently with this assignment,
i.e. so that $v_i<v_j$ whenever $X_{v_i}<X_{v_j}$, then
%
for every vertex $v$ in $G$: 
\begin{description}
\item[(i)] 
$d^r_A(v) \leq 2\frac{d(v)\Delta^{r-1}}{\ln^2\Delta}$,
\item[(ii)] 
$d^r_C(v)\leq 2\frac{d(v)\Delta^{r-1}}{\ln^{3}\Delta}$,
\item[(iii)] 
$\frac{1}{2}\frac{d(v)}{\ln^2\Delta}
\leq d_A(v) \leq 2\frac{d(v)}{\ln^2\Delta}$,
\item[(iv)] 
$\frac{1}{2}\frac{d(v)}{\ln^{3}\Delta}
\leq d_C(v) \leq 2\frac{d(v)}{\ln^{3}\Delta}$,
\item[(v)] 
if $v\in B$, then: $d_-(v)\geq X_v d(v)-\sqrt{X_v d(v)}\ln\Delta$,
\item[(vi)] 
if $v\in B$, then: $d^r_-(v)\leq X_v d(v) \Delta^{r-1}+\sqrt{X_vd(v)\Delta^{r-1}}\ln\Delta$.
\end{description}
\end{lemma}

\begin{pf}
Independently for every $v\in V$ we randomly and uniformly choose a real value $X_v\in [0,1]$
(i.e., we associate with every $v$ an independent random variable $X_v\sim U[0,1]$ having the uniform distribution on $[0,1]$).
With probability one, these values are pairwise distinct for all vertices.
It is also straightforward to note that for every vertex $v\in V$,
${\mathbf E}(d^r_A(v)) \leq \frac{d(v)\Delta^{r-1}}{\ln^2\Delta}$,
${\mathbf E}(d^r_C(v)) \leq \frac{d(v)\Delta^{r-1}}{\ln^{3}\Delta}$,
${\mathbf E}(d_A(v)) = \frac{d(v)}{\ln^2\Delta}$,
${\mathbf E}(d_C(v)) = \frac{d(v)}{\ln^{3}\Delta}$,
${\mathbf E}(d_-(v)) = X_v d(v)$,
${\mathbf E}(d^r_-(v)) \leq X_v d(v) \Delta^{r-1}$.
Then one may prove a concentration of all the corresponding 
random variables using the Chernoff Bound,
which implies that the probability of a contradiction of each of the events \textbf{(i)}--\textbf{(vi)} is bounded from above by $\Delta^{-3r}$.
As each of the $6$ events associated with $v$ is mutually independent of all other such events associated with vertices at distance exceeding 
 $2r$,
analogously as in the previous proof, the thesis is implied by the Lov\'asz Local Lemma, see~\cite{Przybylo_distant_edge_probabil} for details.
\qed
\end{pf}

\section{Proof of Theorem~\ref{przybylo_main_th-sum-proper-probabil}\label{SectionWithPrzybyloProofProbabil}}
Let $r\geq 4$ be a fixed integer and let $G=(V,E)$ be a graph with maximum degree $\Delta$ and
with minimum degree $\delta\geq \ln^8\Delta$.
We shall assume that $\Delta$ is large enough so that all explicit inequalities below hold and $\Delta\geq \Delta'_0$ (from Lemma~\ref{MainSequencingLemma}),
so we shall not specify its value (but assume in particular that $\delta\geq \ln^8\Delta\geq 2$, i.e. there are no isolated edges in $G$).

Let $q$ be the least integer divisible by $3\cdot 2^5 = 96$ such that
\begin{equation}\label{q-bounds}
\frac{\Delta^{r-1}}{\ln\Delta} \leq q < \frac{\Delta^{r-1}}{\ln\Delta}+96,
\end{equation}
and let $Q$ be the least integer divisible by $q$ (thus also by $96$) such that
\begin{equation}\label{Q-bounds}
2\Delta^{r-1}+ \frac{\Delta^{r-1}}{\ln\Delta} \leq Q < 2\Delta^{r-1}+ 2\frac{\Delta^{r-1}}{\ln\Delta} + 96.
\end{equation}



Fix a vertex ordering $v_1,v_2,\ldots,v_n$ of $V$ consistent with Lemma~\ref{MainSequencingLemma} above.
Our goal shall be to show that $\chi'_{\Sigma,r}(G) \leq 2Q+2q$.
For every vertex $v\in A\cup B$ we choose one edge joining it with a vertex in $C$ end denote this edge by $e_v$ -- it exists by \textbf{(iv)} (from  Lemma~\ref{MainSequencingLemma}).
A desired colouring shall be constructed via algorithm developed consistently with the fixed vertex ordering, starting from $v_1$.
Prior launching it we first fix an initial proper edge colouring
$$c_0:E\to\{Q+q-\Delta,Q+q-\Delta+1,\ldots,Q+q\}$$
of $G$, which exists due to the Vizing's Theorem.
Note that this is also a proper edge colouring modulo $q$ (thus also modulo $Q$),
i.e. no two adjacent edges in $G$ have colours congruent modulo $q$.
We shall require this feature within the process of constructing a desired edge colouring from $c_0$, admitting only temporary deviations from this rule or replacing $q$ with $Q$ in the final part of our argument.
While modifying our colouring, by $c(e)$ we shall always mean the contemporary colour of an edge $e$
(hence $d_c(v)$ shall stand for the up-to-date weighted degree of a vertex $v$),
and $d(v)$ shall denote the degree of $v$ in $G$.
In step one of our modifying procedure we shall analyze $v_1$, in step two - $v_2$, and so on.
In general, in step $i$ we shall be modifying only colours of the edges incident with $v_i$ (via rules specified below).
Every vertex $v_i$, the moment it is analyzed (i.e. in step $i$) shall be associated with a $2$-element set, denoted by $S_{v_i}$,
expressing its two admissible sums, 
and belonging to the family (of pairwise disjoint sets):
$$\mathcal{S}=\{\{l,l+Q\}~|~l\in\mathbb{Z} \wedge (l\equiv 0 {\rm ~mod~} {2Q} \vee l\equiv 1 {\rm ~mod~} {2Q} \vee \ldots \vee l\equiv Q-1 {\rm ~mod~} {2Q})\}.$$
Starting from the end of step $i$, we shall require $d_c(v_i)\in S_{v_i}$ till the end of the construction.
The key restriction concerning the choice of 
such set is so that
\begin{description}
\item[($\ast$)] $S_{v_i}$ is disjoint with $S_{v_j}$ for every $j<i$ such that $v_j\in N^r(v_i)$.
\end{description}
This shall be strictly required for all $v_i\in A\cup B$.
%

While modifying colours of the edges, 
we shall obey the following rules.
Suppose a vertex $v$ is being analyzed in a given step. We allow:
\begin{itemize}
  \item[($1^\circ$)] 
  adding $Q$ or subtracting $Q$ (or doing nothing) from the colour of every backward edge of $v$ joining $v$ with a neighbour $u\in A\cup B$
  (so that $d_c(u)\in S_{u}$ afterwards);
\item[($2^\circ$)] 
adding $0$ or $q$ to the colour of every forward edge of $v\in A\cup B$ except $e_v$;
 \item[($3^\circ$)] 
 switching the colour of $e_v$ to any integer in $[Q+q,Q+2q]$ for every $v\in A\cup B$, as long as the edge colouring obtained remains proper modulo $q$. 
\end{itemize}
Note that after introducing such changes 
we shall always have 
\begin{equation}\label{weight_bounds}
q-\Delta\leq c(e)\leq 2Q+2q
\end{equation}
for every $e\in E$ (as desired).
Special rules shall be applied to edges $e$ with both ends in $C$.
These however shall be consistent with~(\ref{weight_bounds}), see details below.
Let us however note here that
by the bounds from (\ref{q-bounds}), (\ref{Q-bounds}) and (\ref{weight_bounds}) above, since $\frac{2Q+2q}{q-\Delta}<5\ln\Delta$, we shall have the following.
\begin{remark}\label{DegreeRemark}
Any $r$-neighbours $u,v$ with
$$d(u) \geq d(v) 5\ln\Delta$$
shall certainly be sum-distinguished in $G$ within 
our construction.
\end{remark}

%

Suppose now we are about to analyze a consecutive
vertex $v\in A$, whose degree we denote by $d$,
and thus far all our rules and requirements have been fulfilled.
Note that
using admissible modifications 
($1^\circ$), ($2^\circ$) and ($3^\circ$)
of the colours of the backward
and forward edges 
of $v$
(since less than $2\Delta$ residues modulo $q$ might be blocked for the colour of $e_v$ due to the required properness of edge colouring modulo $q$),
we may obtain more than
$d (q-2\Delta)$
integer sums at $v$.
At least $d(\frac{q}{3}-2\Delta)$ of these are divisible by $3$.
The set of these (at least) $d(\frac{q}{3}-2\Delta)$ integers contains elements (not necessarily both)
from no less than $d(\frac{q}{6}-\Delta)>2\frac{d\Delta^{r-1}}{\ln^2\Delta}$ pairs from $\mathcal{S}$.
On the other hand, by
\textbf{(i)} (from Lemma~\ref{MainSequencingLemma}), 
$v$ has at most
$\frac{2d\Delta^{r-1}}{\ln^2\Delta}$ backward $r$-neighbours. 
%
We may thus perform admissible alterations of the colours of some of the edges incident with $v$ so that afterwards $d_c(v)$ 
belongs to some pair in $\mathcal{S}$ with elements congruent to $0$ modulo $3$ which is disjoint with 
all $S_u$ associated with backward $r$-neighbours $u$ of $v$.
We set this pair as $S_v$.
We continue in the same manner with all vertices in $A$.

Suppose now that we have reached a vertex $v\in B$ of degree $d$, and thus far all our rules and requirements have been fulfilled.
Similarly as above, admissible modifications ($1^\circ$), ($2^\circ$) and ($3^\circ$) of colours of the edges incident with $v$, due to \textbf{(iv)} and \textbf{(v)}, provide us a list of attainable sums at $v$ of cardinality
(where we in particular additionally use the fact that \textbf{(iv)} implies that $v$ has at least $\frac{d}{2\ln^{3}\Delta} > \frac{Q}{q}$ forward edges.):
\begin{eqnarray}
&&\left(\frac{Q}{q}\left(X_v d-\sqrt{X_v d}\ln\Delta\right)+\left[d-\left(X_v d-\sqrt{X_v d}\ln\Delta\right)\right]\right)\left(q-2\Delta\right)\nonumber\\
&\geq& \left[2\Delta^{r-1}\left(X_v d-\sqrt{X_v d}\ln\Delta\right) + d \frac{\Delta^{r-1}}{\ln\Delta}\right]\left(1-\frac{2\Delta}{q}\right)\nonumber\\
&\geq& \left[2\Delta^{r-1}\left(X_v d-\sqrt{X_v d}\ln\Delta\right) + d \frac{\Delta^{r-1}}{\ln\Delta}\right] - \frac{2\Delta}{q} 2\Delta^{r-1}X_v d -
\frac{2\Delta}{q} d \frac{\Delta^{r-1}}{\ln\Delta}\nonumber\\
&\geq& \left[2\Delta^{r-1}\left(X_v d-\sqrt{X_v d}\ln\Delta\right) + d \frac{\Delta^{r-1}}{\ln\Delta}\right] - 4\Delta \ln\Delta d  -
\frac{1}{2} d \frac{\Delta^{r-1}}{\ln\Delta}\nonumber\\
&\geq& 2\Delta^{r-1}\left(X_v d-\sqrt{X_v d}\ln\Delta\right) + \frac{1}{4} d \frac{\Delta^{r-1}}{\ln\Delta}.\nonumber
\end{eqnarray}
These attainable sums for $v$ contain representatives of at least
$\Delta^{r-1}(X_v d-\sqrt{X_v d}\ln\Delta) + d \frac{\Delta^{r-1}}{8\ln\Delta}$
pairs from $\mathcal{S}$.
On the other hand, \textbf{(vi)} implies that
$$|N^r_-(v)|\leq X_v d \Delta^{r-1}+\sqrt{X_vd\Delta^{r-1}}\ln\Delta,$$
where
$$\Delta^{r-1}(X_v d-\sqrt{X_v d}\ln\Delta) + d \frac{\Delta^{r-1}}{8\ln\Delta} > X_v d \Delta^{r-1}+\sqrt{X_vd\Delta^{r-1}}\ln\Delta.$$
Therefore there is a choice of admissible alterations of the colours of edges incident with $v$ so that afterwards $d_c(v)$
belongs to some set in $\mathcal{S}$ disjoint with $S_u$ for all $u\in N_-^r(v)$.
We perform these alterations and set the corresponding set from $\mathcal{S}$ as $S_v$.
We continue in the same manner with all vertices in $B$.

We are thus left with the analysis of the vertices in $C$.
Let $G'=G[C]$, hence for the maximum degree $\Delta'$ of $G'$ we have $\Delta'\leq \Delta$. Note also that by \textbf{(iv)}, $\delta':=\delta(G')\geq \frac{1}{2}\frac{\delta}{\ln^3\Delta}\geq \frac{1}{2}\ln^5\Delta$. 
Therefore, by Lemma~\ref{SparseSpanningSubgraphInC}, for $\Delta$ sufficiently large,
there exists a spanning subgraph $F'$ of $G'$ with $d_{F'}(v)\leq \frac{d_{G'}(v)}{\ln^3\Delta}$ for every $v\in V$.
Denote the edges of $F'$ by $E'$ 
(hence $F'=(C,E')$),
and note that for every $v\in C$, $d_C(v)-d_{E'}(v)\geq d_{G'}(v)(1-\frac{1}{\ln^3\Delta})\geq 1$ (for $\Delta$ sufficiently large),
hence the edges in $E'':=E(G')\smallsetminus E'=\{e''_1,e''_2,\ldots,e''_m\}$ also induce a spanning subgraph of $G'$.

At this point our edge colouring of $G$ is proper modulo $q$ (hence also modulo $Q$).
We shall now admit a temporary deviation from this rule by setting $c(e)=q$ for every $e\in E'$.
Next we analyze consecutively all edges $e''_1,\ldots,e''_{m''}$ in $E''$
(note that their initial colours, defined by $c_0$, have not been yet altered within our construction, thus all are in the range $[q+Q-\Delta,q+Q]$), and add to a colour of every such subsequent $e''_i=uv$ an integer in $[0,6\Delta]$,
what is consistent with~(\ref{weight_bounds}), so that the obtained sums at $u$ and $v$ are not congruent to $0$ modulo $3$ and so that
the colour of $e''_i$ is not congruent to the colours of its adjacent edges in $G$ modulo $q$.
This is always feasible, as the later requirement blocks at most $2(\Delta-1)$ of at least $2\Delta$ available options in $[0,6\Delta]$ with an adequate residue modulo $3$. After analyzing all edges in $E''$ (inducing a spanning subgraph of $G[C]$), for every vertex $v\in C$ we have $d_c(v)\equiv 1 {\rm ~mod~} {3}$ or $d_c(v)\equiv 2 {\rm ~mod~} {3}$ (contrary to the vertices in $A$).
Now we shall randomly adjust the colours of the edges in $E'$ (which are all set to $q$) to guarantee relatively regular distributions of the sums residues modulo $Q$ in the $r$-neighbourhoods in $C$. In particular we shall show the following.

\begin{lemma}\label{LemmaDistributionInC_new}
We may add to the colour of every edge in $E'$ an integer divisible by $3$ from the set \{0,3,6,\ldots,Q-3\} so that
the obtained edge colouring of $G$ is proper modulo $Q$, 
and
for each vertex $v\in C$ and every integer $t\in [0,Q-1]$ which is not congruent to $0$ modulo $3$,
the number of vertices $u$ in $N_C^r(v)$ with $(5\ln\Delta)^{-1} d(v)\leq d(u) \leq d(v) 5\ln\Delta$
and with $d_c(u)\equiv t {\rm ~mod~} Q$ is upper-bounded by $6000\frac{d(v)}{\ln^{3}\Delta}$. 
\end{lemma}

\begin{pf}
%
%
%
%
%
%
We first partition the set 
$\{0,3,6,\ldots,Q-3\}$ into $32$-element sets of consecutive integers congruent to 0 modulo 3: $L_1,L_2,\ldots,L_{Q/96}$ 
(hence e.g. $L_1=\{0,3,6,\ldots,93\}$).
For every $e\in E'$, as it has less than $2\Delta$ adjacent edges in $G$ (which might block at most $2\Delta$ residues modulo $Q$ for $c(e)$),
i.e. less 
than
$2\Delta$ integers in $[0,Q-1]$ might not be admissible as the additions to the colour of $e$ (equal to $q$ prior to this addition) due to the required properness (modulo $Q$) of the randomly constructed edge colouring.
Thus out of $L_1,L_2,\ldots,L_{Q/96}$,
at least 
$Q/96-2\Delta \geq \frac{\Delta^{r-1}}{48}$ lists (sets)
are entirely available for $e$, where a set $L_i$ is called \emph{entirely available} for $e\in E'$ if neither element of $q+L_i$ is congruent modulo $Q$ to the colour of an edge in $E\smallsetminus E'$ adjacent to $e$ in $G$ (we shall distinguish colours of adjacent edges in $E'$ within our construction below).
Out of these at least $\frac{\Delta^{r-1}}{48}$ entirely available lists for $e$ we randomly and independently for every edge in $E'$ choose one with uniform probability and denote it by $L_e$.
We also temporarily set $c(e)=\min L_e$.

We claim that at the end of such random procedure, with positive probability,
for every $v\in C$ the following event appears:
\begin{description}
\item[$\overline{R_v}$:]
there are at most $31$ edges incident with $v$ (and with both ends in $C$) with a feature that each such edge $e$ is adjacent with an edge $e'$ (with both ends in $C$) such that $L_e=L_{e'}$.
\end{description}
For this goal we shall estimate the probability of the complement of the above for $v\in C$:
\begin{description}
\item[$R_v$:]
there exist $32$ edges incident with $v$ (and with both ends in $C$) with a feature that each such edge $e$ is adjacent with an edge $e'$ (with both ends in $C$) such that $L_e=L_{e'}$.
\end{description}
Fix any $v\in C$ and denote its degree by $d$.
Note first that there are at most ${d \choose 32}\leq {\Delta \choose 32}$ ways of choosing $32$ distinct edges incident with $v$.
Now for a fixed choice of such $32$ edges $B=\{e_1,e_2,\ldots,e_{32}\}$,
each of them is supposed to have an adjacent edge coloured the same (with the same list randomly chosen) as itself,
so for each edge $e_j\in B$ we choose its adjacent edge $e'_j$ which is supposed to have the same colour as $e_j$,
and estimate the probability of $e_1,\ldots,e_{32}$ being witnesses for $R_v$ 
to appear,
by examining all possible configurations of the choices of their correspondents $e'_1,\ldots,e'_{32}$, 
which we divide into $33$ groups
with respect to the number of the edges $e'_j$ belonging to $B$ 
(note that $e'_j$ does not have to be distinct from $e'_l$ for $j\neq l$).
For every $i=0,\ldots,32$ (and fixed $B$), there are at most ${32 \choose i} 31^i (2\Delta)^{32-i}$ choices of edges $e'_1,e'_2,\ldots,e'_{32}$
so that $|\{j:e'_j\in B\}|=i$.
Then for each fixed choice of edges $e'_1,\ldots,e'_{32}$ with this feature, denote $B''=B\cup\{e'_1,e'_2,\ldots,e'_{32}\}$
(hence $32 \leq |B''| \leq 64-i$),
and let us consider an auxiliary graph $H$ with vertex set $B''$ and the set of edges: $\{e_le'_l:l=1,2,\ldots,32\}$.
Note that all its components have order at least $2$.
Fix any subset $B_0\subset B$ of minimal size such that each component of $H$ has at least one vertex in $(B''\smallsetminus B)\cup B_0$,
and note that $|B_0|\leq \lfloor\frac{i}{2}\rfloor$, as there are $32-i$ edges $e_l\in B$  (which are vertices of $H$) adjacent in $H$ with $e'_l\in B''\smallsetminus B$,
while among the remaining at most $i$ edges in $B$ which do not belong to any component including a vertex in $B''\smallsetminus B$
(which induce the remaining components of $H$) it is sufficient to choose at most half to form $B_0$ (one for each of these remaining components of $H$).
Note that edges of $G$ inducing (as vertices o $H$) any component in this auxiliary graph $H$ must have the same colours (lists) chosen to be witnesses for $R_v$ to take place, hence if we fix 
colours (lists) for all edges in $(E'\smallsetminus B)\cup B_0$,
the probability that independent choices for the remaining
at least $32-\lfloor\frac{i}{2}\rfloor$ 
edges in $E'$ 
(from $B\smallsetminus B_0$) shall guarantee $R_v$
is bounded from above by $(\frac{48}{\Delta^{r-1}})^{32-\lfloor\frac{i}{2}\rfloor}$.
By the law of total probability, we thus obtain that: 
\begin{eqnarray}
\mathbf{Pr}(R_v) &\leq& {\Delta \choose 32}\sum_{i=0}^{32}{32 \choose i} 31^i (2\Delta)^{32-i} \left(\frac{48}{\Delta^{r-1}}\right)^{32-\lfloor\frac{i}{2}\rfloor} \nonumber\\
&\leq& 31^{32}\cdot 2^{32}\cdot 48^{32}\Delta^{32} \sum_{i=0}^{32} \Delta^{(32-i)-(32-\lfloor\frac{i}{2}\rfloor)(r-1)} \nonumber\\
&<&10^{48}\cdot 10^{64}\Delta^{32}\cdot 33\Delta^{-16(r-1)} < 10^{114}\Delta^{-4r-12(r-4)}\nonumber\\
&\leq& \frac{10^{114}}{\Delta^{4r}} \label{IneqProbComplRv}
\end{eqnarray}
(for $r\geq 4$).

Now for each vertex $v\in C$ of degree $d$ in $G$ and every integer $t\in [0,Q-1]$ which is not congruent to $0$ modulo $3$,
let $X_{v,t}$ denote (the random variable expressing) the number of vertices $u$ in $N_C^r(v)$ with 
$d_c(u)\in [t-31\cdot 93, t+31\cdot 93] {\rm ~mod~} Q$
(where $c(e)=\min L_e$ for every $e\in E'$)
and $(5\ln\Delta)^{-1} d\leq d(u) \leq d 5\ln\Delta$.
In order to prove the thesis we shall also need to guarantee (with non-zero probability) for every $v\in C$ of degree $d$ (in $G$) and every integer $t\in [0,Q-1]$ which is not congruent to $0$ modulo $3$ the event: 
\begin{description}
\item[$\overline{T_{v,t}}$:] $X_{v,t} \leq 6000 \frac{d}{\ln^{3}\Delta}$.
\end{description}
We thus upper-bound the probability of the complement of this.
As to every edge $e\in E'$ we have assigned the colour being the minimal element $\min L_e$ from the randomly chosen list $L_e$, which may differ by the multiplicity of $96$ between distinct lists, there are at most 
$\lceil (2\cdot 31\cdot 93+1)/96\rceil=61$ distinct values in the interval
$[t-31\cdot 93, t+31\cdot 93]$ the sum at $v$ may possibly attain within our random process.
Therefore for every
$u\in N_C^r(v)$ with $(5\ln\Delta)^{-1} d\leq d(u) \leq d 5\ln\Delta$,
$$\mathbf{Pr}\left(d_c(u) \in [t-31\cdot 93, t+31\cdot 93] {\rm ~mod~} Q\right) \leq 61 \frac{48}{\Delta^{r-1}}$$
(what can be also easily proved by the law of total probability 
via analysis of the possible 
at least $\frac{\Delta^{r-1}}{48}$ choices of lists, hence also additions to the colour, of `the last edge' in $E'$ 
incident with $u$, 
at most $61$ 
of which might assure that
$d_c(u) \in [t-31\cdot 93, t+31\cdot 93] {\rm ~mod~} Q$ regardless of any fixed choices for the remaining edges),
by \textbf{(ii)} 
we thus obtain that
$$\mathbf{E}(X_{v,t})\leq  \frac{61\cdot 48}{\Delta^{r-1}}\frac{2d\Delta^{r-1}}{\ln^{3}\Delta}
= 5856 \frac{d}{\ln^{3}\Delta}.$$
Note also that a change of choice for any edge in $E'$ may influence $X_{v,t}$ by at most $2$.
Moreover, for any $s$, the fact that $X_{v,t}\geq s$ can be certified by the outcomes of at most
$s\cdot \frac{10d}{\ln^{5}\Delta}$ 
trials, i.e. choices committed on the edges in $E'$ incident with some $s$ $r$-neighbours $u$ of $v$ in $C$
with $(5\ln\Delta)^{-1} d\leq d(u) \leq d 5\ln\Delta$, each of which has at most $\frac{2\frac{d 5\ln\Delta}{\ln^3\Delta}}{\ln^3\Delta}=\frac{10d}{\ln^{5}\Delta}$ incident edges in $E'$ by \textbf{(iv)} and Lemma~\ref{SparseSpanningSubgraphInC}.
Thus
by Talagrand's Inequality (and comments below it),
\begin{eqnarray}
\mathbf{Pr}\left(T_{v,t}\right) &\leq&
\mathbf{Pr}\bigg(X_{v,t} > 5856 \frac{d}{\ln^{3}\Delta} + \frac{d}{\ln^{3}\Delta}+ 20\cdot 2\sqrt{\frac{10d}{\ln^{5}\Delta} 5856 \frac{d}{\ln^{3}\Delta}} +64\cdot2^2 \frac{10d}{\ln^{5}\Delta}\bigg)\nonumber\\
&<& 4e^{-\frac{\left(\frac{d}{\ln^{3}\Delta}\right)^2}{8\cdot 2^2 \frac{10d}{\ln^{5}\Delta} \cdot \left(5856 \frac{d}{\ln^{3}\Delta} + \frac{d}{\ln^{3}\Delta}\right)}}<\frac{10^{114}}{\Delta^{4r}}.\label{PrXvtBound}
\end{eqnarray}
%
As any event $T_{v,t}$ 
and $R_v$
is mutually independent of all other events $T_{v',t'}$ 
and $R_{v'}$ with $d(v,v')>2r+1$, i.e., all except at most
$\Delta^{2r+1}\cdot (\frac{2Q}{3}+1)<\Delta^{3r+1}$ such events,
by the Lov\'asz Local Lemma, (\ref{IneqProbComplRv}) and (\ref{PrXvtBound}) we thus obtain that there is a choice of lists (and additions to the colours) of the edges in $E'$ so that none of the events $T_{v,t}$ and $R_v$ holds for any $v\in C$.
This implies among others that each subgraph induced in $G'$ by the edges associated with any fixed list $L_i$ has maximum degree at most $31$.
Thus by Vizing's Theorem we may arbitrarily recolour \emph{properly} each such subgraph, if necessary, using additions 
from its corresponding $L_i$ (where $|L_i|=32$) instead merely the addition $\min L_i$.
Note that then the obtained edge colouring of $G$ is proper modulo $Q$, while colours of some edges could be increased -- each by at most $93$.
As at he same time, every vertex $v$ is by $\overline{R_v}$ incident with at most $31$ edges whose colours could be increased,
by $\overline{T_{v,t}}$ with $v\in C$ and $t\in\{1,2,4,5,7,8,\ldots,Q-2,Q-1\}$ we obtain the thesis.
\qed
\end{pf}


We fix any additions to the colours of the edges in $E'$ 
consistent with the thesis of Lemma~\ref{LemmaDistributionInC_new}.
We shall not alter the colour of any edge with both ends in $C$ anymore, while the remaining ones might be modified by $Q$.
Therefore the edge colouring of $G$ shall remain proper modulo $Q$, while the sums at vertices in $A$ shall remain distinguished from the sums at vertices in $C$, as the first ones are congruent to $0$ modulo $3$, unlike the second ones.
As by \textbf{(iv)} every vertex in $B$ has a neighbour in $C$, we may subtract $Q$ if necessary (or do nothing) from the colour of one such edge for every vertex in $B$ so that the weighted degree for every vertex $v\in B$ is set on the smaller element of its associated two-element list $S_v$.
(This is feasible, as 
prior to these changes, every such edge had its colour between $Q+q-\Delta$ and $Q+2q$,
since it has not been analyzed as a backward edge yet, and therefore 
(\ref{weight_bounds}) shall hold for this edge after any of the described changes).
The thesis of Lemma~\ref{LemmaDistributionInC_new} above obviously still holds afterwards.
The sums at vertices in $B$ shall not be altered anymore.

In the final stage of the construction we shall be subsequently analyzing the vertices in $C$, and modifying colours of the edges joining them with $A$ consistently with ($1^\circ$) in order to dispose of all the remaining sum-conflicts between vertices in $C$ and their $r$-neighbours in $B\cup C$.
This time however we shall admit placing weighted degrees of two $r$-neighbours in the same $2$-element list from $\mathcal{S}$, but in such a way that these weighted degrees are distinct. 
Note that for every consecutive $v\in C$ we have available $d_A(v)+1\geq \frac{d(v)}{2\ln^2\Delta}+1$ (by \textbf{(iii)}) distinct sums, which form an arithmetic progression of difference $Q$, via admissible changes on the edges joining $v$ with $A$. 
These are
all congruent to some $t$ modulo $Q$ (not divisible by $3$) and include at least $\frac{d(v)}{4\ln^2\Delta}$ options which are not fixed as weighted degrees of vertices in $B$, as these are all set to the smaller elements from their associated lists.
So it is sufficient to choose one of such options for $v$ distinct from the contemporary sums at all $r$-neighbours of $v$ in $C$ with
$(5\ln\Delta)^{-1} d(v)\leq d(u) \leq d(v) 5\ln\Delta$ (cf. Remark~\ref{DegreeRemark}) and with weighted degrees congruent to $t$ modulo $Q$.
This is however feasible, as by Lemma~\ref{LemmaDistributionInC_new} above the number of such $r$-neighbours of $v$ equals at most $6000\frac{d(v)}{\ln^{3}\Delta} < \frac{d(v)}{4\ln^2\Delta}$. We choose one of these and perform admissible changes on the edges joining $v$ with $A$
to set it as the sum at $v$.
After analyzing all vertices in $C$, the construction is completed, while the obtained edge colouring $c$ is proper (even modulo $Q$), uses colours in $[q-\Delta,2q+2Q]$ and guarantees sum-distinction between $r$-neighbours in $G$.
\qed

%


\end{document}